\documentclass[draft]{amsart}
\usepackage{mathrsfs}
\usepackage{amssymb}
\usepackage{amsmath}
\usepackage{amsfonts}
\usepackage{amsfonts,enumerate}
\usepackage{pifont}
\usepackage{enumerate}
\usepackage{graphics}
\usepackage{verbatim}
\usepackage{color}
\usepackage{amsthm}

\numberwithin{equation}{section}
\newtheorem{theorem}{Theorem}[section]
\newtheorem{corollary}{Corollary}[section]
\newtheorem{definition}{Definition}[section]

\newtheorem{proposition}{Proposition}[section]
\newtheorem{remark}{Remark}[section]
\newtheorem{example}{Example}[section]

\newcommand{\E}{{\mathcal E}}

\newcommand{\M}{{\mathcal M}}
\newcommand{\N}{{\mathcal N}}

\newcommand{\8}{\infty}
\newcommand{\el}{\ell}

\newcommand{\be}{\begin{eqnarray*}}
\newcommand{\ee}{\end{eqnarray*}}
\newcommand{\beq}{\begin{equation}}
\newcommand{\eeq}{\end{equation}}
\newcommand{\beqn}{\begin{equation*}}
\newcommand{\eeqn}{\end{equation*}}
\newcommand{\bs}{\begin{split}}
\newcommand{\es}{\end{split}}

\begin{document}

\title{Noncommutative martingale inequalities\\ associated with convex functions}


\author[T. N. Bekjan]{Turdebek N. Bekjan}

\address{L. N. Gumilyov Eurasian national University, Astana 010000, Kazakhstan}

\email{bekjant@yahoo.com}


\author[Z. Chen]{Zeqian Chen}

\address{Wuhan Institute of Physics and Mathematics, Chinese
Academy of Sciences, West District 30, Xiao-Hong-Shan, Wuhan 430071, China}

\email{zqchen@wipm.ac.cn}

\date{}

\maketitle

\markboth{T. N. Bekjan, Z. Chen}%
{Noncommutative martingale inequality}

\begin{abstract}
We report recent advances on noncommutative martingale inequalities associated with convex functions. These include noncommutative Burkholder-Gundy inequalities associated with convex functions due to the present authors and Dirksen and Ricard, noncommutative maximal inequalities associated with convex functions due to Os\c{e}kowski and the present authors, and noncommutative Burkholder and Junge-Xu inequalities associated with convex functions due to Randrianantoanina and Lian Wu. Some open problems for noncommutative martingales are also included.
\end{abstract}



\section{Introduction}

The first result of interest on noncommutative martingales was proved by Cuculescu \cite{Cucu1971} in 1971, that is a weak type $(1,1)$ maximal inequality and now known as Cuculescu's construction. But noncommutative martingale theory has received considerable progress since the seminal paper by Pisier and Xu \cite{PX1997} in 1997, thanks to interactions with several fields of mathematics such as operator spaces (e.g. \cite{ER2000, Pisier2003}) and free probability (e.g. \cite{HP2000, VDN1992}). Many classical martingale inequalities have been successfully transferred to the noncommutative setting (cf. e.g. \cite{Xu2007}). Those inequalities of quantum probabilistic nature have, in return, applications to operator spaces (cf. \cite{Junge2005, JP2008, JX2010, LPX2007, PS2002, Xu2006}), quantum stochastic analysis (cf. \cite{JM2010,JM2012,JP2014}) and noncommutative harmonic analysis (cf. \cite{Bek2008,CXY2013,JLeMX2006,JX2007,LMX2010, LMX2011, Mei2007}) et al.

Recall that there are three foundational inequalities in classical martingale theory \cite{Burk1966,Burk1973,BG1970}, which are the Burkholder-Gundy inequality on the square function of martingales, the Doob inequality on the maximal function of martingales, and the Burkholder inequality on the conditioned square function of martingales. These inequalities and their variants characterize martingale Hardy spaces and paly fundamental role in modern analysis, especially in functional analysis, stochastic analysis, harmonic analysis, and mathematical physics, etc. The noncommutative Burkholder-Gundy inequality was established by Pisier and Xu \cite{PX1997} in 1997 as mentioned above, within which they developed a systemical function-analytic method beyond the usual stopping time argument in classical martingale theory. Subsequently, Junge \cite{Junge2002} in 2002 proved the noncommutative Doob inequality via a duality approach. Further, in 2003, Junge and Xu \cite{JX2003} proved the noncommutative Burkholder inequality for characterizing conditioned martingale Hardy spaces, within which, by duality, they presented a dual form of the Burkholder inequality for $1 < p <2.$ This inequality of Junge and Xu is new even in the classical setting, that will be called the Junge-Xu inequality in the sequel for distinguishing it from the Burkholder inequality in $2 \le p < \8.$

This line of investigation was continued by Parcet, Perrin, Randrianantoanina, and the present authors, among others. The weak type $(1,1)$ inequalities for martingale transformations and square and conditioned square functions of noncommutative martingales were proved by Randrianantoanina \cite{Rand2002,Rand2005,Rand2007} between 2002 and 2007. In 2006, Parcet and Randrianantoanina \cite{PR2006} presented Gundy's decomposition for noncommutative martingales and provided alternative proofs of several results mentioned previously. More recently, Perrin, Yin, and the present authors, via a duality approach, have given Davis' and atomic decompositions for noncommutative martingales in $p=1$ in \cite{Perrin2009} and \cite{BCPY2010} respectively. On the other hand, the above three noncommutative martingale inequalities as stated in $L_p$-spaces were generalized to the case of noncommutative symmetric spaces, including noncommutative Lorentz and Orlicz spaces, for details see Dirksen \cite{Dirk2015} and Randrianantoanina and Wu \cite{RandWu2015a} and references therein.

The aim of this paper is to give a survey on noncommutative martingale inequalities associated with convex functions. Recall that classical martingale inequalities associated with convex functions $\Phi$ or $\Phi$-moment martingale inequalities were initiated by Burkholder and Gundy \cite{BG1970}. The above three classical foundational martingale inequalities associated with convex functions were obtained in 1970s' (see \cite{BDG1972,Garsia1973} for details), in which the martingale inequalities in $L_p$-spaces correspond to the case of $\Phi (t) = t^p.$ In 2012, the present authors \cite{BC2012} proved the noncommutative Burkholder-Gundy inequality associated with convex functions, within which a gap in the proof of noncommutative Khintchine's inequality associated with convex functions was fixed later by Dirksen and Ricard \cite{DR2013}. Recently, Os\c{e}kowski and the present authors \cite{BCO2015} proved the noncommtative Doob and ergodic maximal inequalities associated with convex functions. More recently, noncommutative Burkholder and Junge-Xu inequalities associated with convex functions were proved by Randrianantoanina and Wu \cite{RandWu2015b}.

For what follows, in the next section, we will collect notations and definitions for noncommutative $L_p$ and symmetric spaces, convex functions, and noncommutative martingales. In Sections \ref{BurkGundy}, \ref{Doob} and \ref{BurkJungeXu} respectively, we will give detailed descriptions of noncommutative Burkholder-Gundy, Doob, and Burkholder and Junge-Xu inequalities associated with convex functions. Finally, in Section \ref{Open}, we will give some open questions and comments on noncommutative martingale inequalities associated with convex functions.

\section{Preliminaries}\label{pre}

In what follows, by $X \lesssim_{A, B, \ldots} Y$ for two nonnegative (possibly infinite) quantities $X$ and $Y,$ we mean that there exists a constant $C>0$ depending only on $A, B,\ldots$ such that $X \leq C Y,$ and by $X \approx_{A, B} Y$ that $X \lesssim_{A,B} Y$ and $Y \lesssim_{A, B} X.$

\subsection{Symmetric spaces and interpolation with Boyd indices}

First of all, we recall some definitions and notations on symmetric spaces (see \cite{LT1979} for details). We denote by $\mathcal{S}$ the set of all measurable functions $f$ on $(0, \8)$ such that
\be
d_\lambda (f) = | \{ t \in (0, \8): |f(t)| > \lambda \} | < \8
\ee
for some $\lambda <0,$ where $| A|$ is the Lebesgue measure of a measurable subset $A \subset (0, \8).$ For $f \in \mathcal{S},$ we denote by $\mu_t (f)$ the decreasing rearrangement of $f$ in $(0, \8),$ defined by
\be
\mu_t (f) = \inf \{\lambda>0:\; d_\lambda (f) \le t \}
\ee
for all $t >0.$ For $f ,g \in \mathcal{S},$ we say $f$ is submajorized by $g,$ and write $f \prec \prec g,$ if
\be
\int^t_0 \mu_s (f ) d s \le \int^t_0 \mu_s (g) d s
\ee
for all $t > 0.$ A Banach function space $E$ on $(0, \8)$ is called symmetric (resp. fully symmetric) if for $f \in \mathcal{S}$ and $g \in E$ with $\mu (f) \le \mu (g)$ (resp. $f \prec \prec g$), one has $f \in E$ and $\| f \| \le \| g \|.$ A symmetric Banach function space $E$ is said to have the Fatou property if for every increasing net $(f_\alpha)$ in $E$ with $ f_\alpha \ge 0$ and $\sup_\alpha \|f\|_E < \8,$ we have the supremum $f = \sup_\alpha f_\alpha$ exists in $E$ and $\|f_\alpha\|_E \to \| f \|_E.$ The class of (fully) symmetric spaces covers many classical function spaces with the Fatou property, such as Lorentz and Orlicz spaces.

The K\"{o}the dual of a symmetric space $E$ on $(0, \8)$ is the function space defined by setting
\be
E^\times = \Big \{ g \in \mathcal{S}: \; \int | f (t) g (t) | d t < \8,\; \forall f \in E \Big \}
\ee
equipped with the norm $\| g \|_{E^\times} = \sup_{\| f \|_E \le 1} \int | f (t) g (t) | d t.$ Then $E^\times$ is a fully symmetric space having the Fatou property. Note that $E^\times$ is isomorphic to a closed subspace of the Banach dual $E^*$ of $E,$ via the map
\be
g \longrightarrow L_g,\quad L_g (f) = \int f (t) g (t) d t,\; \forall f \in E.
\ee
A symmetric Banach function space $E$ on $(0,\8)$ has the Fatou property if and only if $E = E^{\times \times}$ isometrically. Moreover, a symmetric Banach function space which is separable or has the Fatou property is automatically fully symmetric.

Let $E$ be a symmetric space on $(0, \8).$ For every $0 <s <\8,$ we define $D_s$ in $E$ by
\be
(D_s f) (t) = f (t/s)
\ee
for all $f \in E.$ Then each $D_s$ is a bounded linear operator on $E.$ The lower and upper Boyd indices $p_E$ and $q_E$ of $E$ are defined by
\be\begin{split}
p_E = & \lim_{s \nearrow \8} \frac{\log s}{\log \| D_s \|} = \sup_{s >1} \frac{\log s}{\log \| D_s \|},\\
q_E = & \lim_{s \searrow 0} \frac{\log s}{\log \| D_s \|} = \inf_{0 < s <1} \frac{\log s}{\log \| D_s \|}.
\end{split}\ee
Note that $1 \le p_E \le q_E \le \8,$ and if $E = L^p (0, \8))$ then $p_E = q_E =p$ for all $1 \le p \le \8.$ Note that the
Boyd indices can be expressed as
\be\begin{split}
p_E = & \sup \big \{p>0: \; \exists C_p>0, \; \| D_s \| \le C_p s^{\frac{1}{p}},\;\forall s \ge 1 \big \},\\
q_E = & \inf \big \{q>0: \; \exists C_q>0, \; \| D_s \| \le C_q s^{\frac{1}{q}}, \;\forall s \in (0, 1) \big \}.
\end{split}\ee

Recall that for a given compatible Banach couple $(X, Y)$ in the sense that both are continuously embedded into a Hausdorff topological vector space, a Banach space $Z$ is called an interpolation space if $X \cap Y \subset Z \subset X + Y$ such that whenever a bounded linear operator $T: X + Y \mapsto X + Y$ satisfies $T (X) \subset X$ and $T(Y) \subset Y,$ then $T(Z) \subset Z$ and
\be
\|T\|_Z \le C (\|T \|_X + \|T \|_Y)
\ee
for some constant $C>0.$ In this case, we write $Z \in \mathrm{Int} (X, Y).$ When $C=1,$ $Z$ is called an exact interpolation space. It is known that for a symmetric space $E$ over $(0, \8)$ with $1 \le p < p_E \le q_E < q \le \8,$ we have $E \in \mathrm{Int} (L_p , L_q).$

\subsection{Noncommutative $L_p$ and symmetric spaces}

We use standard notions from theory of noncommutative $L_{p}$-spaces. Our main references are \cite{PX2003} and \cite{Xu2007} (see also \cite{PX2003} for more bibliography). Let $\mathcal{N}$ be a semifinite von Neumann algebra acting on a Hilbert space $\mathbb{H}$ with a normal semifinite faithful trace $\nu.$ A closed densely defined operator $x$ on $\mathbb{H}$ is said to be affiliated with $\N$ if $x$ commutes with every unitary $u$ in the commutant $\N'$ of $\N.$ For a self-adjoint operator $x$ on $\mathbb{H}$ with the spectral decomposition $x = \int^\8_{- \8} t d E_t (x),$ we denote by $\chi_B (x)$ for any Borel subset $B \subset \mathbb{R}$ the the corresponding spectral projection $\int \chi_B (t) d E_t (x).$ An operator $x$ affiliated with $\N$ is called $\nu$-measurable if there exists $s >0$ such
that $\nu [ \chi_{(s,\8)} (|x|)] < \8.$ Let $L_{0}(\mathcal{N})$ denote the topological $*$-algebra of measurable operators with respect to $(\mathcal{N}, \nu).$ The
topology of $L_{0}(\mathcal{N})$ is determined by the convergence in measure. The trace $\nu$ can be extended to the positive cone
$L_{0}^{+}(\mathcal{N})$ of $L_{0}(\mathcal{N}):$
\be
\nu(x)= \int_{0}^{\infty}\lambda\,d\nu(E_{\lambda}(x)),
\ee
where  $x=\int_{0}^{\infty}\lambda\,dE_{\lambda}(x)$ is the spectral decomposition of $x$. Given $0<p<\infty,$ let
\be
L_{p}(\mathcal{N})=\{x\in L_{0}(\mathcal{N}):\;
 \nu(|x|^{p})^{\frac{1}{p}}<\infty\}.
\ee
We define
\be
\|x\|_{p}= \nu(|x|^{p})^{\frac{1}{p}},\quad x\in
L_{p}(\mathcal{N}).
\ee
Then $(L_{p}(\mathcal{N}),\|.\|_{p})$ is a Banach (or quasi-Banach for $p<1$) space. This is the noncommutative $L_{p}$-space associated
with $(\mathcal{N},\nu)$, denoted by $L_{p}(\mathcal{N},\nu)$ or simply by $L_{p}(\mathcal{N}).$ As usual, we set $L_{\infty}(\mathcal{N},\nu)=\mathcal{N}$ equipped with the operator norm.

For $x\in L_{0}(\N)$ we define
\be
\lambda_{s}(x)=\tau(e^{\perp}_s (|x|))\;(s>0)\; \; \text{and}\; \; \mu_t (x) = \inf \{ s>0:\;
\lambda_s (x) \le t \}\; (t >0),
\ee
where $e_s^{\perp} (|x|) = e_{(s,\infty)}(|x|)$ is the spectral projection of $|x|$ associated with the interval $(s,\infty).$ The function $s \mapsto \lambda_{s}(x)$ is called the {\it distribution function} of $x$ and $t \mapsto \mu_{t}(x)$ the {\it generalized singular number} of $x.$ We will denote simply by $\lambda (x)$ and  $\mu(x)$ the functions $s \mapsto \lambda_s (x)$ and $t \mapsto \mu_t (x),$ respectively. It is easy to check that both are decreasing and continuous from the right on $(0,\infty).$ For further information we refer the reader to \cite{FK1986}.

For $0<p<\infty,$ we have the Kolmogorov inequality
\beq\label{eq:kol}
\lambda_{s}(x)\leq\frac{\|x\|_{p}^{p}}{s^{p}},\quad \forall s >0,
\eeq
for any $x\in L_{p}(\N).$ If $x,y$ in $L_0(\N)$, then
\beq\label{eq:tri}
\lambda_{2s}(x+y)\leq \lambda_s(x)+\lambda_s (y),\quad \forall s >0.
\eeq

Throughout, let $E$ be a symmetric space on $(0, \8).$ We define
\be
E(\N, \nu) = \{ x \in L_0 (\N): \; \mu (x) \in E \}
\ee
equipped with the norm $\| x \|_{E(\N)} = \| \mu (x) \|_E.$ The space $E (\N, \nu)$ is a complex Banach space and called noncommutative symmetric space associated with $(\N, \nu)$ relative to $E.$ We refer to \cite{DDP1989, Xu1991, KS2008} for details. We note that if $1 \le p < \8$ and $E= L^p (0, \8),$ then $E (\N, \nu) = L_p (\N, \nu)$ is the above noncommutative $L_p$-space associated with $(\N, \nu).$ It can be shown (cf. \cite[Corollary 2.2]{PX2003}) that if $E$ is a symmetric space over $(0, \8)$ with the Fatou
property such that $1 \le p < p_E \le q_E < q \le \8,$ then for any a semifinite von Neumann algebra $\N,$ $E (\N) \in \mathrm{Int} (L_p (\N), L_q (\N)).$

We define the column space $E (\N, \el^2_c)$ to be the linear space of all sequences $(x_n)_{n \ge 1}$ in $E(\N)$ such that the infinite column matrix
\be
\sum_{n \ge 1} x_n \otimes e_{n, 1} \in E (\N \bar{\otimes} \mathcal{B} (\el^2) )
\ee
where $\mathcal{B} (\el^2)$ is the algebra of all bounded linear operators on $\el^2$ equipped with the usual trace $\mathrm{tr},$ and for $n,m \ge 1,$ $e_{n, m}$ denotes the unit matrix in $\mathcal{B} (\el^2)$ (under the standard basis of $\el^2$) with all entries equal to $0$ but the $(n,m)$ entry being $1.$ Then the linear space
\be
\Big \{ \sum_{n \ge 1} x_n \otimes e_{n, 1}:\; (x_n) \in E (\N, \el^2_c) \Big \}
\ee
is a closed subspace of $E (\N \bar{\otimes} \mathcal{B} (\el^2) ).$ Since the mapping
\be
(x_n) \longmapsto \sum_{n \ge 1} x_n \otimes e_{n, 1}
\ee
from $E (\N, \el^2_c)$ into $E (\N \bar{\otimes} \mathcal{B} (\el^2) )$ is injective, we can equip $E (\N, \el^2_c)$ with the norm
\be
\| (x_n)_{n \ge 1} \|_{E (\N, \el^2_c)} = \Big \| \sum_{n \ge 1} x_n \otimes e_{n, 1} \Big \| = \Big \| \Big ( \sum_{n \ge 1} | x_n|^2 \Big )^\frac{1}{2} \Big \|_{E(\N)}
\ee
so that $E (\N, \el^2_c)$ becomes a Banach space. The corresponding row space $E (\N, \el^2_r)$ is defined to be the linear space of all sequences $(x_n)_{n \ge 1}$ in $E(\N)$ such that $(x^*_n)_{n \ge 1} \in E (\N, \el^2_c),$ equipped with the norm
\be
\| (x_n)_{n \ge 1} \|_{E (\N, \el^2_r)} = \| (x^*_n)_{n \ge 1} \|_{E (\N, \el^2_c)}.
\ee
Then $E (\N, \el^2_r)$ is a Banach space too.

We also need the conditioned versions of the above spaces, which were first introduced by Junge \cite{Junge2002} in $L^p$ cases. To this end, we let $({\mathcal{N}}_{n})_{n \geq 1}$ be an increasing sequence of von Neumann subalgebras of ${\mathcal{N}}$ such that $\cup_{n\geq 1}
{\mathcal{N}}_{n}$ generates ${\mathcal{N}}$ (in the $w^{*}$-topology) and the restriction $\nu_n$ of $\nu$ on $\N_n$ for any $n \ge 1$ is semifinite. For every $n \ge 1,$ we denote by ${\mathcal{E}}_{n}={\mathcal{E}}(\cdot|{\mathcal{N}}_{n})$ the conditional expectation of $(\mathcal{N}, \nu)$ with respect to
$(\mathcal{N}_n, \nu_n)$ with the convention that $\E_0 = \E_1$ below.

For every $n \ge 1$ and $0 < p \le \8,$ we define the conditioned space $L^p_c (\N, \E_n)$ to be the completion of $\N \cap L_2 (\N)$ with respect to the (quasi-)norm
\be
\| x \|_{L^p_c (\E_n)} = \| [ \E_n ( |x|^2)]^\frac{1}{2} \|_p.
\ee
Then, as shown in \cite{Junge2002}, there exists an isometric right $\N_n$-module map $u_{n,p}:\; L^p_c (\N, \E_n) \longmapsto L_p (\N_n, \el^2_c)$ such that
\be
u_{n,p} (x)^* u_{n,q} (y) = \E_n (x^* y) \otimes e_{1,1}
\ee
for all $x \in L^p_c (\N, \E_n)$ and $y \in L^q_c (\N, \E_n)$ with $\frac{1}{p} + \frac{1}{q} \le 1.$

Next, for any finite sequence $(x_n)_{n \ge 1}$ in $L_2 (\N) \cap \N,$ put
\be
U_p[ (x_n)_{n \ge 1} ] = \sum_{n \ge 1} u_{n-1, p} (x_n) \otimes e_{n,1} \in L^p (\N \bar{\otimes} \mathcal{B} (\el^2 (\mathbb{N}^2))).
\ee
Then for any two sequences $(x_n)$ and $(y_n)$ in $\N \cap L_p (\N)$ with $\frac{1}{p} + \frac{1}{q} \le 1,$ one has
\be
U_p[ (x_n)_{n \ge 1} ]^* U_q [ (y_n)_{n \ge 1} ] = \Big ( \sum_{n \ge 1} \E_{n-1} [x^*_n y_n] \Big ) \otimes e_{1,1} \otimes e_{1,1}.
\ee
We define the conditioned column space $L^p_{\mathrm{cond}} (\N, \el^2_c)$ to be the completion of all finite sequences $(x_n)_{n \ge 1}$ in $L_2 (\N) \cap \N$ with respect to the (quasi-)norm
\be
\| (x_n) \|_{L^p_{\mathrm{cond}} (\N, \el^2_c)} = \| U_p [ (x_n)_{n \ge 1} ] \|_{L^p (\N \bar{\otimes} \mathcal{B} (\el^2 (\mathbb{N}^2)))} = \| s^c [(x_n)_{n \ge 1}] \|_p
\ee
where
\be
s^c [(x_n)_{n \ge 1}] = \Big ( \sum_{n \ge 1} \E_{n-1} [ |x_n|^2 ] \Big )^\frac{1}{2}.
\ee
Note that $U_p$ extends to an isometry from $L^p_{\mathrm{cond}} (\N, \el^2_c)$ into $L^p (\N \bar{\otimes} \mathcal{B} (\el^2 (\mathbb{N}^2)))$ and independent of $p$ in the sense of interpolation. We will simply write $U$ for $U_p$ in the sequel.

Since for any finite sequence $(x_n)_{n \ge 1}$ in $L_2 (\N) \cap \N,$
\be
U [ (x_n)_{n \ge 1} ] \in L_2 (\N \bar{\otimes} \mathcal{B} (\el^2 (\mathbb{N}^2))) \cap \N \bar{\otimes} \mathcal{B} (\el^2 (\mathbb{N}^2)),
\ee
we can define the conditioned symmetric column space $E_{\mathrm{cond}} (\N, \el^2_c)$ to be the completion of all finite sequences $(x_n)_{n \ge 1}$ in $L_2 (\N) \cap \N$ with respect to the norm
\be
\| (x_n) \|_{E_{\mathrm{cond}} (\N, \el^2_c)} = \| U [ (x_n)_{n \ge1} ] \|_{E (\N \bar{\otimes} \mathcal{B} (\el^2 (\mathbb{N}^2)))} = \| s^c [(x_n)_{n \ge 1}] \|_{E(\N)}.
\ee
Note that $E_{\mathrm{cond}} (\N, \el^2_c)$ is a Banach space and, $U$ extends to an isometry from $E_{\mathrm{cond}} (\N, \el^2_c)$ into $E (\N \bar{\otimes} \mathcal{B} (\el^2 (\mathbb{N}^2))),$ still denoted by $U.$ Similarly, we can define the conditioned symmetric row space $E_{\mathrm{cond}} (\N, \el^2_r)$ which can also
be viewed as a subspace of $E (\N \bar{\otimes} \mathcal{B} (\el^2 (\mathbb{N}^2)))$ as row vectors.

Note that when $(x_n)_{n \ge 1} \nsubseteqq L_2 (\N) + \N,$ then for any $n \ge 1,$ $|x_n|^2$ may not belong to $L_1 (\N) + \N$ and so $\E_{n-1} (| x_n|^2)$ not be necessarily well-defined. Therefore $s^c [(x_n)_{n \ge 1}] = \sum_{n \ge 1} \E_{n-1} [ |x_n|^2 ]$ is not a well-defined object. To overcome this deficiency, we note that for any finite sequences $(x_n)_{n \ge 1}$ in $L_2 (\N) \cap \N,$
\be
s^c [(x_n)_{n \ge 1}] = | U [ (x_n)_{n \ge 1} ] |.
\ee
Since $U$ extends to an isometry from $E_{\mathrm{cond}} (\N, \el^2_c)$ into $E (\N \bar{\otimes} \mathcal{B} (\el^2 (\mathbb{N}^2))),$ then for any $(x_n)_{n \ge 1} \in E_{\mathrm{cond}} (\N, \el^2_c),$ we can define $s^c [(x_n)_{n \ge 1}]$ as $| U [ (x_n)_{n \ge 1} ] |,$ although $\sum_{n \ge 1} \E_{n-1} [ |x_n|^2 ]$ may not exist as an operator in $L_0 (\N).$ Thus, in general, $\| s^c [(x_n)_{n \ge 1}] \|_{E(\N)}$ is a suggestive notation for every $(x_n)_{n \ge 1} \in E_{\mathrm{cond}} (\N, \el^2_c).$

\subsection{Orlicz functions and noncommutative Orlicz spaces}

Let $\Phi$  be an Orlicz function on $[0,\infty),$ i.e., a continuous increasing and convex function satisfying $\Phi(0)=0$ and $\lim_{t\rightarrow
\infty}\Phi(t)=\infty.$ Recall that $\Phi$ is said to satisfy the $\triangle_2$-condition if there is a constant $C$ such that $\Phi(2t)\leq C\Phi(t)$ for all $t>0.$ In this case, we write $\Phi \in \Delta_2.$ It is easy to check that $\Phi \in \triangle_2$ if and only if for any $a > 0$ there is a constant $C_a>0$ such that $\Phi(a t)\leq C_a \Phi(t)$ for all $t>0.$

For any $x \in L_0 (\N),$ by means of functional calculus applied to the spectral decomposition of $|x|,$ we have
\beq\label{eq:Phispectralintegral}
\nu (\Phi (|x|)) = \int^{\8}_0 \lambda_s (|x|) d \Phi (s) = \int^{\8}_0 \Phi (\mu_t (x)) d t,
\eeq
(see e.g. \cite{FK1986}). Recall that for any $x,y \in L_0 (\N)$ there exist two partial isometries $u,v\in \mathcal{N}$ such that
\beq\label{eq:OperatorModuleTriangleInequa}
| x + y| \le u^* |x| u + v^* |y| v,
\eeq
(cf. \cite{AAP1982}). Then, we have
\be
\nu (\Phi (|\alpha x + (1- \alpha) y|)) \le \alpha \nu (\Phi (|x|)) + (1-\alpha) \nu (\Phi (|y|))
\ee
for any $0 \le \alpha \le 1$ and $x , y \in L_0 (\N).$ In addition, if $\Phi \in \triangle_2,$ then
\be
\nu (\Phi (| x + y|)) \lesssim_\Phi \big [ \nu (\Phi (|x|)) + \nu (\Phi (|y|)) \big ].
\ee

We will work with some standard indices associated to an Orlicz function. Given an Orlicz function $\Phi,$ let
\be M(t, \Phi)= \sup_{s >0} \frac{\Phi (t s)}{\Phi (s)},\quad t >0.
\ee
Define
\be
p_{\Phi} = \lim_{t \searrow 0} \frac{\log M(t, \Phi)}{\log t}, \quad q_{\Phi} = \lim_{t \nearrow \8} \frac{\log M(t, \Phi)}{\log t}.
\ee
These are known as Matuzewska-Orlicz indices of the Orlicz function $\Phi$ (see \cite{M1985, M1989} for more information on Orlicz functions and Orlicz spaces). Note the following properties:
\begin{enumerate}[{\rm 1)}]

\item $1 \le p_{\Phi} \le q_{\Phi} \le \8.$

\item The following characterizations of $p_{\Phi}$ and $q_{\Phi}$ hold
\be p_{\Phi} = \sup \Big \{ p >0:\; \int^t_0 s^{-p} \Phi (s) \frac{d s}{s} = O(t^{- p} \Phi (t)),\; \forall t >0 \Big \};\ee
\be q_{\Phi} = \inf \Big \{ q >0:\; \int^{\8}_t s^{-q} \Phi (s) \frac{d s}{s} = O(t^{- q} \Phi (t)),\; \forall t >0 \Big \}.\ee

\item $\Phi \in \triangle_2$ if and only if $q_{\Phi} < \8,$ or equivalently, $ \sup_{t>0} t \Phi'(t)/\Phi(t)< \8.$ ($\Phi' (t)$ is defined for each $t > 0$ except for a countable set of points in which we take $\Phi'(t)$ as the derivative from the right.)

\end{enumerate}

\begin{example}\label{ex:pleq}\rm
\begin{enumerate}[{\rm i)}]

\item Let $\Phi (t) = t^a \ln (1 + t^b)$ with $a > 1$ and $b >0.$ It is easy to check that $\Phi$ is an Orlicz function and
\be
p_{\Phi} = a\quad \text{and}\quad q_{\Phi} = a + b.
\ee

\item Let $\Phi (t) = t^p (1 + c \sin(p \ln t))$ with $p > 1/(1-2c)$ and $0< c <1/2.$ Then, $\Phi$ is an Orlicz function and
\be
p_{\Phi} = q_{\Phi} = p.
\ee
When $0< c < 1/4,$ $p_{\Phi} = q_{\Phi} =2$ occurs.

\end{enumerate}
\end{example}

For an Orlicz function $\Phi,$ the noncommutative Orlicz space $L_{\Phi}(\mathcal{N})$ is defined as
\be
L_{\Phi}(\mathcal{N}) = \Big \{ x \in L_0 (\N):\; \exists c>0, \nu \Big [ \Phi \Big ( \frac{|x|}{c} \Big ) \Big ] <\infty \Big \}.
\ee
The space $L_{\Phi}(\mathcal{N}),$ equipped with the norm
\be
\|x\|_{\Phi}= \inf \big \{c>0: \;\nu \big [ \Phi({|x|}/{c}) \big ] <1 \big \},
\ee
is a Banach space. If $\Phi(t)=t^p$ with $1 \leq p<\infty$ then $L_\Phi(\mathcal{N})= L_p(\mathcal{N}).$ Note that if $\Phi \in \triangle_2,$ then for $x \in L_0 (\N),$ $x \in L_{\Phi}({\mathcal{N}})$ if and only if $\nu (\Phi (|x|)) < \8.$ Noncommutative Orlicz spaces are fully symmetric spaces of measurable operators as defined in \cite{DDP1989, Xu1991}. Precisely, if $E = L_\Phi (0, \8)$ then $E$ is a (fully) symmetric space on $(0,\8),$ and by \eqref{eq:Phispectralintegral} one has $E (\N) = L_\Phi (\N).$

As described above, we have the Orlicz column and row spaces $L_\Phi (\N, \el^2_c)$ and $L_\Phi (\N, \el^2_r),$ equipped with the norms
\be\begin{split}
\|(x_n)\|_{L_{\Phi}({\mathcal{N}},\el_c^{2})} & = \Big \| \Big ( \sum_n |x_n |^{2} \Big )^{\frac{1}{2}} \Big \|_{\Phi}
\end{split}\ee
and $\| (x_n) \|_{L_{\Phi}({\mathcal{N}},\el_r^{2})} = \| (x^*_n) \|_{L_{\Phi}({\mathcal{N}},\el_c^{2})},$ respectively. Also, we have the conditioned Orlicz column and row spaces $L^\Phi_{\mathrm{cond}} (\N, \el^2_c)$ and $L^\Phi_{\mathrm{cond}} (\N, \el^2_r),$ equipped with the norms respectively
\be\begin{split}
\| (x_n) \|_{L^\Phi_{\mathrm{cond}} (\N, \el^2_c)} & = \| s^c [(x_n)_{n \ge 1}] \|_{\Phi}
\end{split}\ee
and $\| (x_n) \|_{L^\Phi_{\mathrm{cond}} (\N, \el^2_r)} = \| s^c [(x^*_n)_{n \ge 1}] \|_{\Phi}.$

In what follows, unless otherwise specified, we always denote by $\Phi$ an Orlicz function.

\subsection{Noncommutative martingales}

Let ${\mathcal{M}}$ be a semifinite von Neumann algebra with a semifinite normal faithful trace $\tau.$ Let
$({\mathcal{M}}_{n})_{n \geq 1}$ be an increasing sequence of von Neumann subalgebras of ${\mathcal{M}}$ such that $\cup_{n\geq 1}
{\mathcal{M}}_{n}$ generates ${\mathcal{M}}$ (in the $w^{*}$-topology) and for every $n \ge 1,$ the restriction $\tau_n$ of $\tau$ to ${\mathcal{M}}_{n}$ is semifinite, still denoted by $\tau.$ Such a sequence of von Neumann subalgebras of ${\mathcal{M}}$ is called a filtration of $\M.$ Let ${\mathcal{E}}_{n}={\mathcal{E}}(\cdot|{\mathcal{M}}_{n})$ be the conditional expectation of ${\mathcal{M}}$ with respect to
${\mathcal{M}}_{n}.$ Then $\E_n$ extends to a contractive projection from $L_p (\M)$ onto $L_p (\M_n)$ for all $1 \le p \le \8.$ More generally, if $E$ is a symmetric Banach function space on $(0, \8)$ which is an interpolation space of the couple $(L^1 (0, \8), L^\8 (0, \8))$ then $\E_n$ is bounded from $E (\M)$ onto $E(\M_n).$

A noncommutative $E(\M)$-martingale with respect to $({\mathcal{M}}_{n})_{n\geq 1}$ is a sequence $x=(x_{n})_{n\geq 1}$
such that $x_{n} \in E({\mathcal{M}}_{n})$ and
\be
{\mathcal{E}}_n(x_{n+1})=x_n
\ee
for any $n \ge 1.$ Let $\|x\|_{E(\M)}=\sup_{n\geq 1}\|x_{n}\|_{E(\M)}.$ If $\|x\|_{E(\M)} <\infty,$ then $x$ is said to be a bounded $E(\M)$-martingale. If $E = L^1 (0,\8),$ a noncommutative $E(\M)$-martingale is called simply a (noncommutative) martingale.

Let $x$ be a noncommutative martingale. The martingale difference sequence of $x,$ denoted by $dx=(dx_{n})_{n\geq 1},$ is defined as
\be
dx_1=x_1,\quad dx_{n}=x_{n}-x_{n-1},\quad n\geq 2.
\ee
A martingale $x$ is called finite if there exists $N$ such that $d x_n =0$ for all $n \ge N.$ In the sequel, we denote by $x_n = \E_n (x)$ for any $x \in E(\M).$ Note that if $E \in \mathrm{Int} (L^p, L^q)$ with $1 < p \le q < \8$ and satisfies the Fatou property, then any $E(\M)$-bounded martingale $x = (x_n)$ has the form $x_n = \E_n (x_\8)$ for all $n \ge 1,$ where $x_\8 \in E(\M)$ satisfying $\| x_\8 \|_{E(\M)} \approx \| x \|_{E(\M)}$ which is equality if $E$ is an exact interpolation space.

Let $x = (x_n)_{n \ge 1}$ be a $E(\M)$-martingale. Set
\be
S^c (x)= \Big ( \sum_{n = 1}^{\infty} | d x_n |^{2} \Big )^{\frac{1}{2}} \quad \mbox{and} \quad S^r (x)= \Big ( \sum_{n = 1 }^{\infty} | d x_n^* |^{2} \Big )^{\frac{1}{2}}
\ee
where convergences can be taken with respect to the measure topology. These are the column and row versions of the usual square functions for martingales. It should be pointed out that the two operators $S^c (x)$ and $S^r (x)$ may not belong to $E({\mathcal{M}})$ at the same time.

We define $\mathcal{H}^c_E ({\mathcal{M}})$ (resp. $\mathcal{H}^r_E ({\mathcal{M}})$) to be the space of all $E(\M)$-martingales such that $dx \in E ({\mathcal{M}}, \el_c^2 )$ (resp. $dx \in E ({\mathcal{M}}, \el_r^2)$ ), equipped with the norm
\be
\|x\|_{\mathcal{H}^c_E ({\mathcal{M}})}=\|dx\|_{ E ({\mathcal{M}}, \el_c^2)
}
\ee
(resp. $\|x\|_{\mathcal{H}^r_E ({\mathcal{M}})}=\|dx\|_{ E ({\mathcal{M}}, \el_r^2) } $). Since $\mathcal{H}^c_E ({\mathcal{M}})$ and $\mathcal{H}^r_E ({\mathcal{M}})$ embed isometrically into $E (\M \bar{\otimes} \mathcal{B} (\el^2))$ with the closed range and thus they are Banach spaces.

Now, we define the noncommutative symmetric martingale Hardy space $\mathcal{H}_E (\M)$ as follows: If $1 \le p_E \le q_E < 2,$ then
\be
\mathcal{H}_E ({\mathcal{M}}) = \mathcal{H}^c_E ({\mathcal{M}}) + \mathcal{H}^r_E ({\mathcal{M}}),
\ee
equipped with the norm
\be
\|x\|_{\mathcal{H}_E} = \inf \big \{ \|y\|_{ \mathcal{H}^c_E} + \|z\|_{\mathcal{H}^r_E}:\;
x=y+z,\; y \in \mathcal{H}^c_E ({\mathcal{M}}),\; z \in \mathcal{H}^r_E ({\mathcal{M}}) \big \}.
\ee
If $2 \leq p_E \le q_E < \8,$
\be
\mathcal{H}_E ({\mathcal{M}}) = \mathcal{H}^c_E ({\mathcal{M}}) \cap \mathcal{H}^r_E ({\mathcal{M}}),
\ee
equipped with the norm
\be
\|x\|_{\mathcal{H}_E} = \|x\|_{\mathcal{H}^c_E} + \|x\|_{\mathcal{H}^r_E}.
\ee

If $E= L^p (0, \8)$ with $1 < p < \8,$ then the noncommutative Burkholder-Gundy inequality states that
\be
\| x \|_p \approx \| x \|_{\mathcal{H}_p}
\ee
which was proved by Pisier and Xu \cite{PX1997} in 1997.

Next, we consider noncommutative conditioned Hardy spaces developed in \cite{JX2003}. For any finite martingale $x = (x_n)_{n\ge1}$ in $L_2(M)$, we set
\be
s^c (x)= s^c [(d x_n)_{n \ge 1}] = \Big ( \sum_{n \ge 1 }\E_{n-1}[ |dx_n |^{2}] \Big )^{\frac{1}{2}}
\ee
and $s^r (x) = s^c (x^*) = s^c [(d x^*_n)_{n \ge 1}],$ respectively. These are called the column and row conditioned square functions. We denote by $\mathrm{h}^c_E ({\mathcal{M}})$ and $\mathrm{h}^r_E ({\mathcal{M}})$ the corresponding completion spaces of all finite martingale $x = (x_n)_{n\ge1}$ in $L_2(M)$ under the norms
\be\begin{split}
\|x\|_{\mathrm{h}^c_E ({\mathcal{M}})} & = \| (d x_n)_{n \ge 1} \|_{E_{\mathrm{cond}} (\M, \el^2_c)} = \|s^c (x)\|_{E(\M)},\\
\|x\|_{\mathrm{h}^r_{E}({\mathcal{M}})} & = \| (d x_n)_{n \ge 1} \|_{E_{\mathrm{cond}} (\M, \el^2_r)} = \|s^r (x)\|_{E(\M)},
\end{split}\ee
respectively. We call $\mathrm{h}^c_E ({\mathcal{M}})$ and $\mathrm{h}^r_E ({\mathcal{M}})$ the conditioned column and row martingale Hardy spaces.

Also, we define the diagonal martingale Hardy space $\mathbf{h}^d_E (\M)$ to be the space of all martingales whose martingale difference sequences belong to $E (\M \bar{\otimes} \el^\8)$ equipped with the norm $\| x \|_{\mathbf{h}^d_E} = \| (d x_n) \|_{E (\M \bar{\otimes} \el^\8)}.$

Now we are ready to define the conditioned martingale Hardy space $\mathbf{h}_E.$ If $1 \le p_E \le q_E < 2,$ then
\be
\mathbf{h}_E ({\mathcal{M}}) = \mathbf{h}^d_E (\M) + \mathbf{h}^c_E ({\mathcal{M}}) + \mathbf{h}^r_E ({\mathcal{M}}),
\ee
equipped with the norm
\be
\|x\|_{\mathbf{h}_E} = \inf \big \{ \|w\|_{\mathbf{h}^d_E} + \|y\|_{\mathbf{h}^c_E} + \|z\|_{\mathbf{h}^r_E}:\;
x=w + y+z,\; w \in \mathbf{h}^d_E (\M), y \in \mathbf{h}^c_E ({\mathcal{M}}),\; z \in \mathbf{h}^r_E ({\mathcal{M}}) \big \}.
\ee
If $2 \leq p_E \le q_E < \8,$
\be
\mathcal{H}_E ({\mathcal{M}}) = \mathbf{h}^d_E (\M) \cap \mathbf{h}^c_E ({\mathcal{M}}) \cap \mathbf{h}^r_E ({\mathcal{M}}),
\ee
equipped with the norm
\be
\|x\|_{\mathbf{h}_E} = \|x\|_{\mathbf{h}^d_E} +  \|x\|_{\mathbf{h}^c_E} + \|x\|_{\mathbf{h}^r_E}.
\ee

If $E= L^p (0, \8)$ with $1 < p < \8,$ then the noncommutative Burkholder and Junge-Xu inequality states that
\be
\| x \|_p \approx \| x \|_{\mathbf{h}_p}
\ee
which was proved by Junge and Xu in \cite{JX2003}. Note that the case $1 < p <2$ is new even in the classical setting and now known as the Junge-Xu inequality, as compared with the original Burkholder inequality in $2 \le p < \8.$

\section{Noncommutative Burkholder-Gundy inequalities}\label{BurkGundy}

The proof of noncommutative Burkholder-Gundy inequalities associated with convex functions is based noncommutative Marcinkiewicz interpolation theorem and Khintchine inequalities for Rademacher's random variables associated with convex functions, which are of independent interest.

Let $\mathcal{N}_1$ (resp. $\mathcal{N}_2$) be a semifinite von Neumann algebra on a Hilbert space $\mathbb{H}_1$ (resp. $\mathbb{H}_2$)
with a normal semifinite faithful  trace $\nu_1$ (resp. $\nu_2$). Recall that a map $T:L_{0}(\mathcal{N}_1)\rightarrow L_{0}(\mathcal{N}_2)$ is said to
be sublinear if for any  operators $x,y\in L_{0}(\mathcal{N}_1),$ there exist two partial isometrics $u,v\in \mathcal{N}_2$ such that
$$
|T(x+y)|\leq u^{*}|Tx|u+v^{*}|Ty|v,\quad |T(\alpha x)|\leq
|\alpha||Tx|, \;\forall\alpha\in \mathbb{C}.
$$
Note that for any $x,y \in L_0 (\N)$ there exist two partial isometrics $u,v\in \mathcal{N}$ such that
\beq\label{eq:OperatorModuleInequa}
| x + y| \le u^* |x| u + v^* |y| v,
\eeq
(see \cite{AAP1982}) and then a linear operator is sublinear. We recall that a sublinear operator $T: L_{0}(\mathcal{N}_1)\rightarrow L_{0}(\mathcal{N}_2)$ is of {\it weak type} $(p,q)$ with $1 \le p \le q \le \8,$ if there is a constant $C>0,$ so that for every $x \in L_p (\N_1)$
\beq\label{eq:WLp}
\lambda_{\alpha}(|Tx|) \leq \Big ( \frac{C\|x\|_p}{\alpha} \Big )^q,\quad \forall \alpha > 0.
\eeq
If $q = \8,$ it means that $\| T x \|_q \le C \| x \|_p.$

The classical Marcinkiewicz interpolation theorem has been extended to include Orlicz spaces as interpolation classes by A.Zygmund, A.P.Calder\'{o}n, S.Koizumi, I.B.Simonenko, W.Riordan, H.P.Heinig and A.Torchinsky (for references see \cite{M1989}). The following result is a noncommutative Marcinkiewicz interpolation theorem associated with convex functions, which seems new even in the classical setting, at least we have no concrete reference.

\begin{theorem}\label{th:InterBG} {\rm (\cite[Theorem 2.1]{BC2012})}\;
Let $\mathcal{N}_1$ (resp. $\mathcal{N}_2$) be a semifinite von Neumann algebra on a Hilbert space $\mathbb{H}_1$ (resp. $\mathbb{H}_2$) with a normal semifinite faithful  trace $\nu_1$ (resp. $\nu_2$).
Suppose $1 \le p_{0}<p_{1} \leq \infty.$ Let $T:L_{0}(\mathcal{N}_1) \rightarrow
L_{0}(\mathcal{N}_2)$ be a sublinear operator and simultaneously of weak types $(p_i, p_i)$ for $i=0$ and $i=1.$ If $\Phi$ is an Orlicz function with $p_{0}<p_{\Phi}\le q_{\Phi}<p_{1},$ then there exists a constant $C$ depending only on $p_0, \; p_1$ and $\Phi,$ such that
\begin{equation}\label{strongphi}
\nu_2 (\Phi(|Tx|))\leq C \nu_1 ( \Phi(|x|)),
\end{equation}
for all $x \in L_\Phi(\mathcal{N}_1).$
\end{theorem}

\begin{remark}\rm
If $T$ is of strong type $(p, p),$ i.e., there exists a constant $C>0$ such that $\| T x \|_p \le C \| x \|_p$ for any $x \in L_p (\N),$ then by the Kolmogorov inequality \eqref{eq:kol} we have \be \lambda_{\alpha}(|Tx|)\leq \alpha^{-p} \| Tx \|^p_p \le C^p \alpha^{-p} \|x \|^p_p,\ee
that is, $T$ is of weak type $(p, p).$ Consequently, if $T$ is simultaneously of strong types $(p_i, p_i)$ for $i=0$ and $i=1,$ then the conclusion of Theorem \ref{th:InterBG} holds.
\end{remark}

Recall that a sequence of random variables $\{\varepsilon_i\}$ on a probability space $(\Omega, P)$ is called a Rademacher's sequence if it is independent and satisfies
\be
P (\varepsilon_i =1) = P (\varepsilon_i = -1) = \frac{1}{2}
\ee
for all $i.$ The following is noncommutative Khintchine's inequalities for Rademacher's sequences associated with convex functions.

\begin{theorem}\label{th:khin} {\rm \cite[Theorem 4.1]{BC2012}}\;
Let $\M$ be a semifinite von Neumann algebra with a semifinite normal faithful trace $\tau.$ Let $\Phi$ be an Orlicz function and $\{\varepsilon_i\}$ a Rademacher's sequence on a probability space $(\Omega, P).$
\begin{enumerate}[\rm (1)]

\item If $1<p_{\Phi} \le q_{\Phi}<2,$ then for any finite sequence $\{x_k\}$ in $L_{\Phi}(\M),$
\begin{equation}\label{eq:khin1}
\begin{split}
\int_{\Omega} \tau \Big ( \Phi \Big [ & \Big | \sum_{k=0}^{n}x_{k}\varepsilon_{k} \Big | \Big ] \Big ) d P \\
& \approx_\Phi \inf \Big \{ \tau \Big ( \Phi \Big [ \Big ( \sum_{k=0}^{n}
|y_{k}|^{2} \Big )^{\frac{1}{2}} \Big ] \Big )+ \tau \Big ( \Phi \Big [ \Big ( \sum_{k=0}^{n}
|z_{k}^{*}|^{2} \Big )^{\frac{1}{2} } \Big ] \Big ) \Big \},
\end{split}
\end{equation}
where the infimun runs over all decomposition $x_{k}=y_{k}+z_{k}$ with $y_{k}$ and $z_{k}$ in $L_{\Phi}({\M}).$

\item If $2<p_{\Phi} \le q_{\Phi}<\infty,$ then for any finite sequence $\{x_k\}$ in $L_{\Phi}(\M),$
\begin{equation}\label{eq:khin2}
\begin{split}
\int_{\Omega} \tau \Big ( \Phi \Big [ & \Big | \sum_{k=0}^{n}x_{k}\varepsilon_{k} \Big | \Big ] \Big ) d P\\
& \approx_\Phi \tau \Big ( \Phi \Big [ \Big ( \sum_{k=0}^{n}
|x_{k}|^{2} \Big )^{\frac{1}{2} } \Big ] \Big ) + \tau \Big ( \Phi \Big [ \Big ( \sum_{k=0}^{n}
|x_{k}^{*}|^{2} \Big )^{\frac{1}{2} } \Big ] \Big ).
\end{split}
\end{equation}
\end{enumerate}
\end{theorem}

\begin{remark}\rm
\begin{enumerate}[ 1)]

\item Note that Khintchine's inequality is valid for $L_1$-norm in both commutative and noncommutative settings (cf., \cite{LPP1991}). We could conjecture that the right condition in Theorem \ref{th:khin} (1) should be $q_{\Phi} <2$ without the additional restriction condition $1<p_{\Phi}.$ However, our argument seems to be inefficient in this case. We need new ideas to approach it.

\item We note that there is a gap in the proof of the inequality $``\le"$ in \eqref{eq:khin2} in \cite{BC2012}, as pointed out to the present authors by Q. Xu. This was resolved subsequently by Dirksen and Ricard in \cite{DR2013}.

\end{enumerate}
\end{remark}

Now, we are in a position to state the noncommutative Burkholder-Gundy martingale inequalities associated with convex functions.

\begin{theorem}\label{th:BG}{\rm \cite[Theorem 5.1]{BC2012}}\;
Let $\M$ be a semifinite von Neumann algebra with a semifinite normal faithful trace $\tau$ and $({\mathcal{M}}_{n})_{n\geq 1}$ a filtration
of ${\mathcal{M}}.$ Let $\Phi$ be an Orlicz function and $x = ( x_n )_{n\geq 1}$ a noncommutative $L_{\Phi}$-martingale with respect to $({\mathcal{M}}_{n})_{n\geq 1}.$
\begin{enumerate}[\rm (1)]

\item If $1<p_{\Phi} \le q_{\Phi}<2,$ then
\begin {equation}\label{eq:BG1}
\begin{split}
\tau \big ( \Phi [ |x|] \big ) \approx_\Phi \inf \Big \{ \tau \Big ( \Phi \Big [ \Big ( \sum_{n= 1}^{\infty} |d y_n |^{2} \Big )^{ \frac{1}{2}} \Big ] \Big ) + \tau \Big ( \Phi \Big [ \Big ( \sum_{n= 1}^{\infty} |d z_n^{*}|^{2} \Big )^{ \frac{1}{2}} \Big ] \Big ) \Big \}
\end{split}
\end{equation}
where the infimum runs over all decomposition $x_n = y_n + z_n$ with $(y_n)_{n \ge 1}$ in $\mathcal{H}_{C}^{\Phi}({\mathcal{M}})$ and $(z_n)_{ne \ge 1}$ in
$\mathcal{H}_{R}^{\Phi}({\mathcal{M}}).$

\item If $2 <p_{\Phi} \le q_{\Phi}<\infty,$ then
\begin {equation}\label{eq:BG2}
\tau \big ( \Phi [ |x|] \big ) \approx_\Phi \tau \Big ( \Phi \Big [ \Big ( \sum_{n = 1}^{\infty} |d x_n |^{2} \Big )^{ \frac{1}{2}} \Big ] \Big ) + \tau \Big ( \Phi \Big [ \Big ( \sum_{n = 1}^{\infty} |d x_n^{*}|^{2} \Big )^{ \frac{1}{2}} \Big ] \Big ).
\end{equation}
\end{enumerate}
\end{theorem}

\begin{remark}\rm
Theorem \ref{th:BG} was sharpened recently by Y. Jiao, F. Sukochev, D. Zanin, and D. Zhou \cite{JSZZ2015}, as informed to the second-named author by Prof. Sukochev in 2015 summer seminar in Wuhan University.
\end{remark}

\section{Noncommutative maximal inequalities}\label{Doob}

The noncommutative Doob and ergodic maximal inequalities associated with convex functions are based on a noncommutative Marcinkiewicz interpolation theorem for convex functions of maximal operators. To this end, we need to introduce the definition of convex functions of maximal operators in the noncommutative setting.

\subsection{Noncommutative maximal functions in $L_p$-spaces}

Given $1 \le p < \8,$ $L_p(\M; \el^{\8})$ is defined as the space of all sequences $(x_n)_{n \ge 1}$ in $L_p (\M)$ for which there exist $a, b \in L_{2p}(\M)$ and a bounded sequence $(y_n)_{n \ge 1}$ in $\M$ such that $x_n = a y_n b$ for all $n \ge 1.$ For such a sequence, set
\beq\label{eq:p-MaxNorm}
\left \| ( x_n )_{n \ge 1} \right \|_{L_p (\M, \el^{\8})} : = \inf \big \{ \| a \|_{2 p} \sup_n\| y_n \|_{\8} \| b \|_{2 p} \big \},
\eeq
where the infimum runs over all possible factorizations of $(x_n)_{n \ge 1}$ as above. This is a norm and $L_p (\M; \el^{\8})$ is a Banach space. These spaces were first introduced by Pisier \cite{Pisier1998} in the case when $\M$ is hyperfinite and by Junge \cite{Junge2002} in the general case, and studied extensively by Junge and Xu \cite{JX2007}.

As in \cite{JX2007}, we usually write
\be
\big \| {\sup_n}^+ x_n \big \|_p = \| ( x_n )_{n \ge 1} \|_{L_p (\M, \el^{\8})}.
\ee
We warn the reader that this suggestive notation should be treated with care. It is used for
possibly nonpositive operators and
\be
\big \| {\sup_n}^+ x_n \big \|_p \neq \big \| {\sup_n}^+ | x_n | \big \|_p
\ee
in general. However it has an intuitive description in the positive case, as observed in \cite[p.329]{JX2007}: A positive sequence
$(x_n)_{n \ge1}$ of $L_p (\M )$ belongs to $L_p (\M; \el^{\8})$ if and only if there exists a positive $a \in L_p (\M)$ such
that $x_n \le a$ for any $n \ge 1$ and in this case,
\beq\label{eq:PositiveSequence}
\big \| {\sup_n}^+ x_n \big \|_p = \inf \big \{\| a \|_p :\; a \in L_p (\M),\; x_n \le a, \; \forall n \ge 1 \big \}.
\eeq
In particular, it was proved in \cite{JX2007} that the spaces $L_p (\M; \el^{\8})$ for all $1 \le p \le \8$ form interpolation scales with respect to complex interpolation. However, this result is no longer true for the real interpolation. This is one of the
difficulties one will encounter for dealing with Marcinkiewicz type interpolation theorem on maximal operators in the noncommutative setting.

\begin{definition}\label{df:QuasiOperator}
Let $1 \le p_0 < p_1 \le \8.$ Let $S = (S_n)_{n \ge 1}$ be a sequence of maps from $L^+_{p_0}(\mathcal{M}) + L^+_{p_1}(\mathcal{M}) \mapsto L^+_{0}(\mathcal{M}).$
\begin{enumerate}[\rm (1)]

\item $S$ is said to be subadditive, if for any $n \ge 1,$
\be
S_n (x + y) \le S_n (x) + S_n (y),\quad \forall x, y \in L^+_{p_0}(\mathcal{M}) + L^+_{p_1}(\mathcal{M}).
\ee

\item $S$ is said to be of weak type $(p, p)$ ($p_0 \le p < p_1$) if there is a positive constant $C$ such that for any $x \in L^+_p (\M)$ and any $\lambda >0$ there exists a projection $e \in \M$ such that
\be
\tau (e^{\perp} ) \le \left ( \frac{C \| x \|_p}{\lambda} \right )^p \quad \text{and} \quad e S_n (x) e \le \lambda,\; \forall n \ge 1.
\ee

\item $S$ is said to be of type $(p, p)$ ($p_0 \le p \le p_1$) if there is a positive constant $C$ such that for any $x \in L^+_p (\M)$ there exists $a \in L^+_p (\M)$ satisfying
\be
\| a \|_p \le C \| x \|_p \quad \text{and} \quad S_n (x) \le a,\; \forall n \ge 1.
\ee
In other words, $S$ is of type $(p, p)$ if and only if $\| S(x) \|_{L_p (\M; \el^{\8})} \le C \| x \|_p$ for all $x \in L^+_p (\M).$

\end{enumerate}
\end{definition}

This definition of subadditive operators in the noncommutative setting is due to Junge and Xu \cite{JX2007}, who proved a noncommutative analogue of the classical Marcinkiewicz interpolation theorem as follows.

\begin{theorem}\label{th:InterJX2007}{\rm (cf. \cite[Theorem 3.1]{JX2007})}\;
Let $1 \le p_0 < p_1 \le \8.$ Let $S = (S_n)_{n \ge 1}$ be a sequence of maps from $L^+_{p_0}(\mathcal{M}) + L^+_{p_1}(\mathcal{M}) \mapsto L^+_{0}(\mathcal{M}).$ Assume that $S$ is subadditive. If $S$ is of weak type $(p_0, p_0)$ with constant $C_0$ and of type $(p_1, p_1)$ with constant $C_1,$ then for any $p_0 < p < p_1,$ $S$ is of type $(p, p)$ with constant $C_p$ satisfying
\be
C_p \le C C^{1- \theta}_0 C^{\theta}_1 \left ( \frac{1}{p_0} - \frac{1}{p} \right )^{- 2}
\ee
where $\theta$ is determined by $1 /p = (1- \theta) /p_0 + \theta /p_1$ and $C$ is an absolute constant.
\end{theorem}

\subsection{Interpolation for convex functions of noncommutative maximal operators}

To state our results, we need to define a convex function of maximal operators in the noncommutative setting as follows.

\begin{definition}\label{df:PhiMaxOper}{\rm \cite[Definition 3.2]{BCO2015}}\;
Let $\Phi$ be an Orlicz function. Let $(x_n)$ be a sequence in $L_{\Phi} (\M).$ We define $\tau [ \Phi ( {\sup_n}^+ x_n ) ]$ by
\beq\label{eq:PhiMaxSequ}
\tau \Big [ \Phi \big ( {\sup_n}^+ x_n \big ) \Big ] : = \inf \left \{ \frac{1}{2} \Big ( \tau \big [ \Phi \big ( |a|^2 \big ) \big ] + \tau \big [ \Phi \big ( |b|^2 \big ) \big ]\Big ) \sup_n \| y_n \|_{\8} \right \}
\eeq
where the infimum is taken over all decompositions $x_n = a y_n b$ for $a, b \in L_0 (\M)$ and $(y_n) \subset L_{\8} (\M)$ with $|a|^2, |b|^2 \in L_{\Phi} (\M),$ and $\| y_n \|_{\8} \le 1$ for all $n.$
\end{definition}

To understand $\tau \big [ \Phi \big ( {\sup_n}^+ x_n \big ) \big ],$ let us consider a positive sequence $x = (x_n)$ in $L_{\Phi} (\M).$ We then note that
\beq\label{eq:PhiMaxSequPosiGe}
\tau \Big [ \Phi \big ( {\sup_n}^+ x_n \big ) \Big ] \le \tau \big [ \Phi \big ( a \big ) \big ],
\eeq
if $a \in L_{\Phi}^+ (\M)$ such that $x_n \le a$ for all $n.$ Indeed, for every $n$ there exists a contraction $u_n$ such that $x_n^{\frac{1}{2}} = u_n a^{\frac{1}{2}}$ and hence $x_n = a^{\frac{1}{2}} u^*_n u_n a^{\frac{1}{2}}.$ This concludes \eqref{eq:PhiMaxSequPosiGe}. Moreover, the converse to \eqref{eq:PhiMaxSequPosiGe} also holds true provided $\Phi \in \triangle_2.$

\begin{proposition}\label{prop:BasicPhiMax}
Let $\Phi$ be an Orlicz function satisfying the $\triangle_2$-condition.
\begin{enumerate}[\rm (1)]

\item If $x = (x_n)$ is a positive sequence in $L_{\Phi} (\M),$ then
\be
\tau \Big [ \Phi \big ( {\sup_n}^+ x_n \big ) \Big ] \approx \inf \Big \{ \tau \big [ \Phi \big ( a \big ) \big ]:\; a \in L_{\Phi}^+ (\M)\; \text{such that}\; x_n \le a, \forall n \ge 1 \Big \}.
\ee

\item For any two sequences $x =(x_n), y = (y_n)$ in $L_{\Phi} (\M)$ one has
\be
\tau \Big [ \Phi \big ( {\sup_n}^+ (x_n + y_n) \big ) \Big ] \lesssim \tau \Big [ \Phi \big ( {\sup_n}^+ x_n \big ) \Big ] + \tau \Big [ \Phi \big ( {\sup_n}^+ y_n \big ) \Big ].
\ee


\end{enumerate}
\end{proposition}

\begin{remark}\label{rk:MaxOrliczSpace}\rm
For a sequences $x =(x_n)$ in $L_{\Phi} (\M),$ set
\be
\big \| {\sup_n}^+ x_n \big \|_{\Phi} : = \; \inf \Big \{ \lambda > 0:\; \tau \Big [ \Phi \Big ( {\sup_n}^+ \frac{x_n}{\lambda} \Big ) \Big ] \le 1 \Big \}.
\ee
One can check that $\| {\sup_n}^+ x_n \|_{\Phi}$ is a norm in $x = (x_n).$ Define
\be
L_{\Phi} (\M; \el^{\8}): = \; \Big \{(x_n) \subset L_{\Phi} (\M):\; \tau \Big [ \Phi \Big ( {\sup_n}^+ \frac{x_n}{\lambda} \Big ) \Big ] < \8\; \text{for some}\; \lambda >0 \Big \},
\ee
equipped with $\| (x_n) \|_{L_{\Phi} (\M; \el^{\8})} = \| {\sup_n}^+ x_n \|_{\Phi}.$ Then $L_{\Phi} (\M; \el^{\8})$ is a Banach space. For $1 \le p < \8,$ if $\Phi (t) = t^p$ then $L_{\Phi} (\M; \el^{\8}) = L_p (\M; \el^{\8})$ with equivalent norms. The details are left to the interested readers.
\end{remark}

We are ready to state the noncommutative Marcinkiewicz type interpolation theorem for a convex function of maximal operators.

\begin{theorem}\label{th:InterDoob}{\rm \cite[Theorem 3.2]{BCO2015}}\;
Let $S = (S_n)_{n \ge 0}$ be a sequence of maps from $L^+_1(\mathcal{M}) + L^+_{\8}(\mathcal{M}) \mapsto L^+_{0}(\mathcal{M}).$ Let $1 \le p < \8.$ Assume that $S$ is subadditive, and order preserving in the sense that for all $n \ge 1,$ $S_n(x) \le S_n(y)$ whenever $x \le y$ in $L^+_0 (\M).$ If $S$ is simultaneously of weak type $(p, p)$ with constant $C_p$ and of type $(\8, \8)$ with constant $C_{\8},$ then for an Orlicz function $\Phi$ with $p <p_{\Phi}\le q_{\Phi}< \8,$ there exists a positive constant $C$ depending only on $C_p, C_{\8}, p_{\Phi}$ and $q_{\Phi},$ such that
\begin{equation}\label{eq:PhiMaxOper}
\tau \Big [ \Phi \big ( {\sup_n}^+ S_n( x) \big ) \Big ]\leq C \tau \big [ \Phi( x) \big ]
\end{equation}
for all $x \in L^+_{\Phi}(\mathcal{M}).$
\end{theorem}

\begin{remark}\label{re:InterDoob}
\begin{enumerate}[\rm (1)]\rm

\item The classical Marcinkiewicz interpolation theorem has been extended to include Orlicz spaces as interpolation classes by A. Zygmund, A. P. Calder\'{o}n {\it et al.} (cf. e.g. \cite{M1989}). The noncommutative analogue of this associated with a convex function was recently obtained in \cite{BC2012}. Theorem \ref{th:InterDoob} can be considered as a noncommutative Marcinkiewicz type interpolation theorem for a convex function of maximal operators, which seems new even in the classical setting.

\item Should Theorem \ref{th:InterDoob} be true whenever $S$ is simultaneously of weak type $(p_0, p_0)$ and of type $(p_1, p_1)$ and $\Phi$ an Orlicz function such that $1 \le p_0 <p_{\Phi}\le q_{\Phi}< p_1 \le \8$ (i.e., the case $p_1<\8$ is included). This was raised as an open question in the preliminary version of \cite{BCO2015}, and affirmatively answered by Dirksen \cite{Dirk2012b} recently through extending the argument of \cite{BCO2015} to the case of restricted weak type inequalities for noncommutative maximal operators. However, Theorem \ref{th:InterDoob} is sufficient for proving noncommutative maximal inequalities associated with convex functions (see Theorems \ref{th:PhiDoob} and \ref{th:NcMaxErgodi} below).

\end{enumerate}
\end{remark}

The argument presented in the proof of Theorem \ref{th:InterDoob} can be used to obtain the corresponding interpolation theorem for maximal operators on noncommutative symmetric spaces, which is different from that of \cite{Dirk2015} and slightly simpler.

\begin{theorem}\label{th:InterSymspace}{\rm \cite[Theorem 3.3]{BCO2015}}\;
Let $S = (S_n)_{n \ge 0}$ be a sequence of maps from $L^+_1(\mathcal{M}) + L^+_{\8}(\mathcal{M}) \mapsto L^+_{0}(\mathcal{M}).$ Assume that $S$ is subadditive and order preserving. Let $1 \le p < \8.$ Let $E$ be a rearrangement invariant space such that $p_E > p.$ If $S$ is simultaneously of weak type $(p, p)$ with constant $C_p$ and of type $(\8, \8)$ with constant $C_{\8},$ then there exists a positive constant $C_E$ depending only on $C_p, C_{\8}, p$ and $p_E,$ such that for any $x \in E^+ (\M, \tau)$ there exists $a \in E^+ (\M, \tau)$ satisfying
\begin{equation}\label{eq:InterSymspace}
\| a \|_{E (\M, \tau)} \le C_E \| x \|_{E(\M, \tau)} \quad \text{and}\quad S_n (x) \le a, \quad \forall n \ge 0.
\end{equation}
\end{theorem}

\begin{remark}\label{re:InterSymspace}\rm
With the help of Theorem \ref{th:InterSymspace}, the associated maximal inequalities on noncommutative symmetric spaces are in order, including Doob's inequality, Dunford-Schwartz and Stein maximal ergodic inequalities, as well as the corresponding pointwise convergence theorems (see \cite{JX2007} for detailed information).
\end{remark}

\subsection{Noncommutative Doob and ergodic maximal inequalities associated with convex function}

The first main result is the following noncommutative Doob inequality associated with a convex function, generalizing Junge's noncommutative Doob inequality in $L_p$ in \cite{Junge2002}.

\begin{theorem}\label{th:PhiDoob}{\rm \cite[Theorem 4.1]{BCO2015}}\;
Let $\M$ be a semifinite von Neumann algebra with a semifinite normal faithful trace $\tau,$ equipped with a filtration $( \M_n )_{n\ge 0}$ of von Neumann subalgebras of $\M.$ Let $\Phi$ be an Orlicz function and $x = ( x_n )$ be a noncommutative $L_{\Phi}$-martingale with respect to $({\mathcal{M}}_{n}).$ If $1 < p_{\Phi} \le q_{\Phi}< \8,$ then
\beq\label{eq:PhiDoob}
\tau \Big [ \Phi \big ( {\sup_n}^+ x_n \big ) \Big ] \thickapprox_\Phi \tau \big [ \Phi( |x| ) \big ].
\eeq
\end{theorem}

\begin{remark}\label{re:PhiMaxSpace}\rm
Let $\Phi$ be an Orlicz function. We define the Hardy-Orlicz maximal space of noncommutative martingales as
\be
\mathcal{H}^{\mathrm{max}}_{\Phi} (\M): = \Big \{ x \in L_{\Phi} (\M):\; \| x \|_{\mathcal{H}^{\mathrm{max}}_{\Phi}} = \big \| {\sup_n}^+ \mathcal{E}_n (x) \big \|_{\Phi} < \8 \Big \}.
\ee
(See \cite[Sect. 4]{JX2005} for the case $\Phi (t) = t^p$.) Then, Theorem \ref{th:PhiDoob} implies that $\mathcal{H}^{\mathrm{max}}_{\Phi} (\M) = L_{\Phi} (\M)$ with equivalent norms provided $1 < p_{\Phi} \le q_{\Phi}<\8.$
\end{remark}

As a consequence of Theorem \ref{th:PhiDoob}, we obtain the following noncommutative Burkholder-Davis-Gundy inequality associated with a convex function.

\begin{corollary}\label{cor:NcBDG}{\rm \cite[Corollary 4.1]{BCO2015}}\;
Let $\M$ be a semifinite von Neumann algebra with a semifinite normal faithful trace $\tau,$ equipped with a filtration $( \M_n )_{n\ge 0}$ of von Neumann subalgebras of $\M.$ Let $\Phi$ be an Orlicz function, and let $x=( x_{n} )_{n\geq 0}$ be a noncommutative $L_{\Phi}$-martingale with respect to $( \M_n )_{n \ge 0}.$ If $1<p_{\Phi} \le q_{\Phi}<2,$ then
\beq\label{eq:NcPhiBDGle2}
\begin{split}
\tau \Big ( \Phi \Big [ {\sup_n}^+ x_n \Big ] \Big ) \approx_\Phi \inf \bigg \{ \tau \Big ( \Phi \Big [ \Big ( \sum_{k= 0 }^{\infty} |dy_{k}|^{2}
\Big )^{ \frac{1}{2}} \Big ] \Big ) + \tau \Big ( \Phi \Big [ \Big ( \sum_{k= 0 }^{\infty} |dz_{k}^{*}|^{2} \Big )^{ \frac{1}{2}} \Big ] \Big ) \bigg \},
\end{split}
\eeq
where the infimum runs over all decomposition $x_n = y_n + z_n$ with $y_n$ in $\mathcal{H}^C_{\Phi}({\mathcal{M}})$ and $z_n$ in $\mathcal{H}^R_{\Phi}({\mathcal{M}});$ and if $2 <p_{\Phi} \le q_{\Phi}<\infty,$ then
\beq\label{eq:NcPhiBDGge2}
\begin{split}
\tau \Big ( \Phi  \Big [ {\sup_n}^+ x_n \Big ] \Big ) \approx_\Phi \tau \Big ( \Phi \Big [ \Big ( \sum_{k= 0 }^{\infty} |dx_{k}|^{2} \Big )^{ \frac{1}{2}} \Big ] \Big ) + \tau \Big ( \Phi \Big [ \Big ( \sum_{k= 0 }^{\infty} |dx_{k}^{*}|^{2} \Big )^{ \frac{1}{2}} \Big ] \Big ).
\end{split}\eeq
\end{corollary}

\begin{remark}\label{rk:NcBDG}\rm
The classical case of Corollary \ref{cor:NcBDG} was originally proved by Burkholder, Davis, and Gundy in \cite{BDG1972} (see also \cite{Burk1973}). Note that, the classical case holds even if $p_\Phi =1$ (e.g. \cite{Davis1970}). However, the noncommutative case is surprisingly different. Indeed, it was shown in \cite[Corollary 14]{JX2005} that $\mathcal{H}_1 \not= \mathcal{H}^{\mathrm{max}}_1.$ This implies that \eqref{eq:NcPhiBDGle2} does not hold when $ \Phi (t) = t$ for which $p_\Phi =1.$
\end{remark}

Now we turn to noncommutative maximal ergodic inequalities associated with convex functions. To state our results, we need some more notation.

Let $\M$ be a semifinite von Neumann algebra with a normal semifinite faithful trace $\tau,$ and let $L_p(\M)$ be the associated noncommutative $L_p$-space. Consider a linear map $T:\; \M \mapsto \M$ which may satisfy the
following conditions:
\begin{enumerate}[{\rm (I)}]

\item $T$ is a contraction on $\M,$ that is, $\| T x \| \le \|x\|$ for all $x \in \M.$

\item $T$ is positive, i.e., $Tx \ge 0$ if $x \ge 0.$

\item $\tau \circ T \le \tau,$ that is, $\tau (T x) \le \tau (x)$ for all $x \in L_1 (\M) \cap \M_+.$

\item $T$ is symmetric relative to $\tau,$ i.e., $\tau ( (Ty)^* x ) = \tau (y^* T x )$ for all $x, y \in L_2 (\M) \cap \M.$

\end{enumerate}
Under conditions (I)-(III), $T$ naturally extends to a contraction on $L_p ( \M )$ for every $1 \le p < \8.$
The extension will be still denoted by $T.$

\begin{theorem}\label{th:NcMaxErgodi}{\rm \cite[Theorem 4.2]{BCO2015}}\;
Let $\Phi$ be an Orlicz function with $1 < p_{\Phi} \le q_{\Phi} < \8.$ If $T: \M \mapsto \M$ is a linear map satisfying $(I) - (III),$ then
\beq\label{eq:NcPhiDS}
\tau \Big ( \Phi \big [ {\sup_n}^+ M_n (x) \big ] \Big ) \lesssim_\Phi \tau ( \Phi [ |x| ] ),\quad \forall x \in L_{\Phi} (\M),
\eeq
where $M_n : = \frac{1}{n+1} \sum^n_{k=0} T^k$ for any $n \ge 1.$ If, in addition, $T$ satisfies (IV), then
\beq\label{eq:NcPhiStein}
\tau \Big ( \Phi \big [ {\sup_n}^+ T^n (x) \big ] \Big ) \lesssim_\Phi \tau ( \Phi [ |x| ] ),\quad \forall x \in L_{\Phi} (\M).
\eeq
\end{theorem}

The inequalities \eqref{eq:NcPhiDS} and \eqref{eq:NcPhiStein} are the noncommutative forms of the classical Dunford-Schwartz and Stein maximal ergodic inequality for a convex function of positive and symmetric positive contractions. These generalize the noncommutative Dunford-Schwartz and Stein maximal ergodic inequalities of Junge and Xu \cite{JX2007} in the $L^p$ case to the case of convex functions. The proofs of \eqref{eq:NcPhiDS} and \eqref{eq:NcPhiStein} are again based on Theorem \ref{th:InterDoob}.

\subsection{Noncommutative weak type maximal inequalities}


All the results continue to hold if we replace the noncommutative maximal operator $\tau [ \Phi ({\sup_n}^+ x_n) ]$ associated with a convex function by a certain weak maximal operator. The required modifications are not difficult and left to the interested reader. However, for
the sake of convenience, we write the corresponding definitions and results. We refer to \cite{BCLJ2011} for noncommutative weak Orlicz spaces and for the terminology used here.

Let $\Phi$ be an Orlicz function. For $x \in L^w_{\Phi} (\M),$ we set
\be
\| x \|_{\Phi, \8} = \sup_{t >0} t \Phi \big [ \mu_t (x) \big ].
\ee
When $\Phi (t) = t^p,$ $\| x \|_{\Phi, \8}$ is just the usual weak $L_p$-norm $\| x \|_{p, \8}.$

The following is the definition of a weak type maximal operator associated with a convex function.

\begin{definition}\label{df:PhiMaxOperW}{\rm \cite[Definition 5.1]{BCO2015}}\;
Let $(x_n)$ be a sequence in $L^w_{\Phi} (\M).$ We define $\| {\sup_n}^+ x_n \|_{\Phi, \8}$ by
\beq\label{eq:PhiMaxSequ}
\big \| {\sup_n}^+ x_n \big \|_{\Phi, \8} : = \inf \left \{ \frac{1}{2} \Big ( \big \| |a|^2 \big \|_{\Phi, \8} +  \big \| |b|^2 \big \|_{\Phi, \8} \Big ) \sup_n \| y_n \|_{\8} \right \}
\eeq
where the infimum is taken over all decompositions $x_n = a y_n b$ for $a, b \in L_0 (\M)$ and $(y_n) \subset L_{\8} (\M)$ with $|a|^2, |b|^2 \in L^w_{\Phi} (\M),$ and $\| y_n \|_{\8} \le 1$ for all $n.$
\end{definition}

We have the noncommutative Marcinkiewicz type interpolation theorem for the weak type maximal operator associated with a convex function as follows. To this end, recall that
\be
a_{\Phi} = \inf_{t>0} \frac{t \Phi'(t)}{\Phi (t)}\quad \text{and} \quad b_{\Phi} = \sup_{t>0} \frac{t \Phi'(t)}{\Phi (t)}.
\ee
Note that $1 \le a_\Phi \le p_\Phi \le q_\Phi \le b_\Phi \le \8,$ but they do not coincide in general (see \cite{M1985, M1989} for details).

\begin{theorem}\label{th:InterDoobW}{\rm \cite[Theorem 5.1]{BCO2015}}\;
Suppose $1 \le p_{0}<p_{1} \leq \infty.$ Let $S = (S_n)_{n \ge 0}$ be a sequence of maps from $L^+_{p_0}(\mathcal{M}) + L^+_{p_1}(\mathcal{M}) \mapsto L^+_{0}(\mathcal{M}).$ Assume that $S$ is subadditive. If $S$ is of weak type $(p_0, p_0)$ with constant $C_0$ and of type $(p_1, p_1)$ with constant $C_1,$ then for an Orlicz function $\Phi$ with $p_{0}<a_{\Phi}\le b_{\Phi}<p_{1},$ there exists a positive constant $C$ depending only on $p_0, \; p_1, C_0, C_1$ and $\Phi,$ such that
\begin{equation}\label{eq:PhiMaxOperW}
\big \| {\sup_n}^+ S_n( x) \big \|_{\Phi, \8} \leq C \left \| x \right \|_{\Phi, \8},
\end{equation}
for all $x \in L^w_{\Phi}(\mathcal{M})_+.$
\end{theorem}

The following is a noncommutative Doob weak type inequality associated with a convex function.

\begin{theorem}\label{th:PhiDoobW}{\rm \cite[Theorem 5.2]{BCO2015}}\;
Let $\M$ be a semifinite von Neumann algebra with a semifinite normal faithful trace $\tau,$ equipped with a filtration $({\mathcal{M}}_{n})$ of von Neumann subalgebras of ${\mathcal{M}}.$ Let $\Phi$ be an Orlicz function and let $x = ( x_n )$ be a noncommutative $L^w_{\Phi}$-martingale with respect to $({\mathcal{M}}_{n}).$ If $1 < a_{\Phi} \le b_{\Phi}< \8,$ then
\beq\label{eq:PhiDoobW}
\big \| {\sup_n}^+ x_n \big \|_{\Phi, \8} \approx_\Phi \| x \|_{\Phi, \8}.
\eeq
\end{theorem}

The weak type analogue of Theorem \ref{th:NcMaxErgodi} concerning maximal ergodic inequalities associated with a convex function is similar and omitted.

\section{Noncommutative Burkholder and Junge-Xu inequalities}\label{BurkJungeXu}

Burkholder's inequality \cite{Burk1973} is a fundamental result in classical martingale theory, stating that for every $p >2,$ the $L_p$-space is equal to the intersection space of corresponding conditioned and diagonal Hardy spaces. In a remarkable paper \cite{JX2003}, Junge and Xu obtained the noncommutative analogue of this inequality for characterizing noncommutative conditioned martingale Hardy spaces $\mathbf{h}_p.$ By duality, they proved a corresponding inequality yet for $p<2,$ showing that noncommutative $L_p$-space coincides with the sum space of $\mathbf{h}_p$ and the diagonal martingale Hardy space $\mathbf{h}^d_p$ of noncommutative martingales if $1 < p<2.$ This inequality of Junge and Xu is new even for classical (commutative) martingales. Noncommutative Burkholder inequality was extended recently by Dirksen \cite{Dirk2015} to the case of symmetric spaces, including noncommutative Orlicz spaces. More recently, Randrianantoanina and Wu \cite{RandWu2015a} proved the Junge-Xu inequality in noncommutative symmetric spaces.

In this section, we will describe the very recent result of Randrianantoanina and Wu \cite{RandWu2015b} about noncommutative Burkholder and Junge-Xu inequalities associated with convex functions. To sate their result, we need to introduce some notations. Let $\Phi$ be an Orlicz function. As shown in Section \ref{pre}.2, the conditioned square function $s^c (x)$ is in general not well defined for a $L_\Phi$-martingale and so does $\tau [\Phi (s^c(x))].$ To overcome this obstacle, note that the mapping $D: x \longmapsto (d x_n)_{n \ge 1}$ from finite $L_2$-martingales into $L_\Phi (\M, \el^2_c)$ extends to an isometry from $\mathbf{h}^\Phi_c (\M)$ into $L_\Phi (\M, \el^2_c).$ Recall that $U$ is a isometry from $L_\Phi (\M, \el^2_c)$ into $L_\Phi (\M \bar{\otimes} \el^2 (\mathbb{N}^2))$ (cf. Section \ref{pre}.2). Then, we can define
\be
s^c (x) = |U D (x) |
\ee
for all $x \in \mathbf{h}^\Phi_c (\M),$ since it holds for all finite $L_2$-martingales $x.$ Therefore, we define
\be
\tau [\Phi (s^c(x))] = \tau \otimes \mathrm{Tr} [\Phi (|U D (x) |)]
\ee
for every $x \in \mathbf{h}^\Phi_c (\M),$ where $\mathrm{Tr}$ denotes the usual trace on $ \mathcal{B} (\el^2 (\mathbb{N}^2)).$

The following is the noncommutative Burkholder and Junge-Xu inequalities associated with convex functions.

\begin{theorem}\label{thm:Burk}{\rm \cite{RandWu2015b}}\;
Let $\M$ be a semifinite von Neumann algebra with a semifinite normal faithful trace $\tau,$ equipped with a filtration $(\M_n)_{n \ge 1}$ of it. Let $\Phi$ be an Orlicz function with the Boyd indexes $p_\Phi$ and $q_\Phi.$ Then for every noncommutative $L_{\Phi}$-martingale $x=(x_{n})_{n \geq 1},$ if $2 <p_{\Phi} \le q_{\Phi}< \8,$
\beq\label{eq:Burkconvfunct}
\begin{split}
\tau \big [ \Phi ( |x|) \big ] \approx_\Phi \sum_{n \ge 1} \tau [ \Phi ( |dx_n|) ] + \tau [ \Phi ( s^c (x))] + \tau [ \Phi (s^r (x)) ];
\end{split}\eeq
and if $1 < p_{\Phi} \le q_{\Phi}< 2,$
\beq\label{eq:JungeXuconvfunct}
\begin{split}
\tau \big [ \Phi ( |x|) \big ] \approx_\Phi \inf \Big \{ \sum_{n \ge 1} \tau [ \Phi ( |dw_n|) ] + \tau [ \Phi ( s^c (y) ) ] + \tau [ \Phi ( s^r (z)) ] \Big \},
\end{split}\eeq
where the infimum is taken over all the decompositions $x = w + y + z$ with $w \in \mathbf{h}^d_\Phi,$ $y \in \mathbf{h}^c_\Phi,$ and $z \in \mathbf{h}^r_\Phi.$
\end{theorem}

The proof of Randrianantoanina and Wu in \cite{RandWu2015b} was primarily motivated by an observation that
\be
\nu (\Phi (|x|)) = \int^{\8}_0 \Phi (\mu_t (x)) d t \thickapprox \int^\8_0 \Phi [ t^{-1} K (t; L_1 (\M), \M ) ] d t,
\ee
where $K (t; L_1 (\M), \M )$ is the $K$-functional from interpolation theory. They first proved the Junge-Xu inequality \eqref{eq:JungeXuconvfunct} by heavily employing results from interpolation theory and then deduce the Burkholder inequality \eqref{eq:Burkconvfunct} from \eqref{eq:JungeXuconvfunct} via a duality argument adapted to convex functions in a nontrivial way.

\section{Concluding remarks and open problems}\label{Open}

In this section, we collect some remarks and open problems to our knowledge.

\begin{remark}\label{rk:EndBG}\rm
The first question concerns with the endpoint case of the noncommutative Burkholder-Gundy inequality \eqref{eq:BG1}. Precisely, we would conjecture that the following is correct: For an Orlicz function $\Phi$ satisfying the $\triangle_2$-condition, there is a constant $C>0$ such that if $x = ( x_n )_{n \ge 1}$ is a  noncommutative $L_1$-bounded martingale, then there exist two martingales $y$ and $z$ with $x = y +z$ satisfying
\be
\Phi \big [ \mu_t \big ( S^c (y) \big ) \big ] + \Phi \big [ \mu_t \big ( S^r (z) \big ) \big ] \le \frac{C}{t} \| \Phi (|x|) \|_1
\ee
for any $t >0.$ The endpoint case of the noncommutative Junge-Xu inequality associated with convex functions seems open yet.
\end{remark}

\begin{remark}\label{rk:EndBurk}\rm
The second question is raised by Randrianantoanina and Wu \cite{RandWu2015b} concerning with the noncommutative Burkholder inequality \eqref{eq:Burkconvfunct}. That is, does the Burkholder inequality $`` \le "$ in \eqref{eq:Burkconvfunct} hold under the weaker condition $1 < p_\Phi \le q_\Phi < \8$ ?

Note that Dirksen and Ricard proved in \cite{DR2013} that if $\Phi$ is an Orlicz function with $1 < p_\Phi \le q_\Phi < \8,$ then
\be
\tau [\Phi (|x|)] \lesssim_\Phi \tau \big [ \Phi ( S^c (x) ) \big ] + \tau \big [ \Phi ( S^c (x^*) ) \big ]
\ee
for all $x \in L_\Phi (\M).$ This improves the inequality $``\le"$ in \eqref{eq:BG2}.
\end{remark}

\begin{remark}\label{rk:MaxBurk}\rm
The third question concerns with the mixture of the noncommutative Doob and Burkholder inequality associated with convex functions. Precisely, if $\Phi$ is an Orlicz function with $1 < p_\Phi \le q_\Phi < \8,$ do we have
\be
\tau [\Phi (|x|)] \lesssim_\Phi \tau \big [ \Phi ( {\sup}^+_n d x_n ) \big ] + \tau \big [ \Phi ( s^c (x) ) \big ] + \tau \big [ \Phi ( s^c (x^*) ) \big ]
\ee
for all $x \in L_\Phi (\M)$ ? This problem is closely related with \cite[Problem 6.3]{RandWu2015b}.
\end{remark}

\begin{remark}\label{rk:TypeIIIOrlicz}\rm
Recently, Labuschagne in \cite{Lab2013} has presented a noncommutative Orlicz space theory for type III von Neumann algebras. As pointed out in \cite{RandWu2015b}, it is interesting to generalize noncommutative martingale inequalities associated with convex functions previously stated for semifinite von Neumann algebras to the case of type III von Neumann algebras.
\end{remark}

\begin{remark}\label{rk:atom}\rm
The most open problem in noncommutative martingale theory should be the atomic decomposition problem raised in \cite{BCPY2010}, i.e., whether or not $\mathbf{h}^c_p = \mathbf{h}^{c, \mathrm{at}}_p$ for $0 < p <1$ ? We refer to \cite[Appendix A]{BCPY2010} and \cite{HM2012} for details. We have known that whenever $\mathbf{h}^c_p$ admits a $(p,2)_c$-atomic decomposition, then it must possess a $(p, \8)_c$-atomic decomposition (cf. \cite{Chen}).
\end{remark}

\subsection*{Acknowledgement} The authors are grateful to Prof. Xu for helpful suggestion and encouragement. T. N. Bekjan is partially supported by NSFC grant No.11371304. Z. Chen is in part supported by NSFC grant No.11171338 and No.11431011.


\begin{thebibliography}{0}

\bibitem{AAP1982} C. A. Akemann, J. Anderson, G. K. Pedersen,
Triangle inequalities in operator algebras,
{\it Linear Multilinear Algebra} {\bf 11} (1982), 167-178.



\bibitem{Bek2008} T. N. Bekjan,
Noncommutative maximal ergodic theorems for positive contractions,
{\it J. Funct. Anal.} {\bf 254} (2008), 2401-2418.

\bibitem{BC2012} T. N. Bekjan, Z. Chen,
Interpolation and $\Phi$-moment inequalities of noncommutative martingales,
{\it Probab. Theory Relat. Fields} {\bf 152} (2012), 179-206.

\bibitem{BCLJ2011} T. N. Bekjan, Z. Chen, P. Liu, Y. Jiao,
Noncommutative weak Orlicz spaces and martingale inequalities,
{\it Studia Math.} {\bf 204}(3) (2011), 195-212.

\bibitem{BCPY2010} T. N. Bekjan, Z. Chen, M. Perrin, Z. Yin,
Atomic decomposition and interpolation for Hardy spaces of noncommutative martingales,
{\it J. Funct. Anal.} {\bf 258} (2010), 2483-2505.

\bibitem{BCO2015} T. N. Bekjan, Z. Chen, A. Os\c{e}kowski,
Noncommutative maximal inequalities associated with convex functions,
to appear in {\it Trans. Amer. Math. Soc.}




\bibitem{Burk1966} D. L. Burkholder.
Martingale transforms.
{\it Ann. math. Statist.} {\bf 37} (1966), 1494-1504.

\bibitem{Burk1973} D. L. Burkholder,
Distribution function inequalities for martingales,
{\it Ann. Probab.} {\bf 1}(1) (1973), 19-42.


\bibitem{BDG1972} D. L. Burkholder, B. Davis, R. Gundy,
Integral inequalities for convex functions of operators on martingales, in
{\it Proc. 6th Berkley Symp.} {\bf II},223-240, 1972.

\bibitem{BG1970} D. L. Burkholder, R. Gundy.
Extrapolation and interpolation of quasi-linear operators on martingales.
{\it Acta Math.} {\bf 124} (1970), 249-304.

\bibitem{Chen} Z. Chen, in preparation.

\bibitem{CXY2013} Z. Chen, Q. Xu, Z. Yin,
Harmonic analysis on quantum tori,
{\it Commun. Math. Phys.} {\bf 322} (2013), 755-805.

\bibitem{Cucu1971} I. Cuculescu,
Martingales on von Neumann algebras,
{\it J. Multivariate. Anal.} {\bf 1} (1971), 17-27.

\bibitem{Davis1970} B. J. Davis,
On the integrability of the martingale square function,
{\it Israel J. Math.} {\bf 8} (1970), 187-190.


\bibitem{Dirk2015} S. Dirksen,
Noncommutative Boyd interpolation theorems,
{\it Trans. Amer. Math. Soc.} {\bf 367} (6) (2015), 4079-4110.

\bibitem{Dirk2012b} S. Dirksen,
Weak-type interpolation for noncommutative maximal operators,
arXiv:1212.5168, 2012.

\bibitem{DR2013} S. Dirksen, E. Ricard,
Some remarks on noncommutative Khintchine inequalities,
{\it Bull. London Math. Soc.} {\bf 45} (2013), 618-624.

\bibitem {DDP1989} P. G. Dodds, T. K. Dodds, B. de Pagter,
Noncommutative Banach function spaces,
{\it Math. Z. } {\bf 201} (1989), 583-587.

\bibitem {DDP1992} P. G. Dodds, T. K. Dodds, B. de Pagter,
Fully symmetric operator spaces,
{\it Integ. Equ. Oper. Theory} {\bf 15} (1992), 942-972.

\bibitem{FK1986} T. Fack, H. Kosaki,
Generalized $s$-numbers of $ \tau $-measure operators,
{\it Pacific J. Math. } {\bf 123} (1986), 269-300.

\bibitem{ER2000} E. Effros, Z. J. Ruan,
{\it Operator Spaces,}
Oxford University Press, Oxford, 2000.

\bibitem{Garsia1973} A. M. Garsia.
On a convex function inequality for martingales.
{\it Ann. Probab.} {\bf 1} (1973), 171-174.


\bibitem{HP2000} F. Hiai, D. Petz,
{\it The Semicircle Law, Free Random Variables and Entropy,}
Amer. Math. Soc., Providence, RI, 2000.

\bibitem{HM2012} G. Hong, T. Mei,
John-Nirenberg inequalities and atomic decomposition for noncommutative martingales,
{\it J. Funct. Anal.} {\bf 263} (2012), 1064-1097.



\bibitem{JSZZ2015} Y. Jiao, F. Sukochev, D. Zanin, D. Zhou,
Johnson-Schechtman inequalities for noncommutative martingales, in preprint.

\bibitem{Junge2002} M. Junge,
Doob's inequality for non-commutative martingales,
{\it J. Reine Angew. Math.} {\bf 549} (2002), 149-190.

\bibitem{Junge2005} M. Junge,
Embedding of the operator space OH and the logarithmic `little Grothendieck inequality',
{\it Invent. Math.} {\bf 161} (2005), 225-286.

\bibitem{JLeMX2006} M. Junge, C. Le Merdy, Q. Xu,
$\mathrm{H}^{\8}$ functional calculus and square functions on noncommutative $L_p$-spaces,
{\it Ast\'{e}risque} {\bf 305} (2006), vi + 138.

\bibitem{JM2010} M. Junge, T. Mei,
Noncommutative Riesz transforms - a probabilistic approach,
{\it Amer. J. Math.} {\bf 132} (2010), 611-681.

\bibitem{JM2012} M. Junge, T. Mei,
BMO spaces associated with semigroups of operators,
{\it Math. Ann.} {\bf 352} (2012), 691-743.

\bibitem{JMP2014} M. Junge, T. Mei, J. Parcet,
Smooth Fourier multipliers on group von Neumann algebras,
{\it Geometric Funct. Anal.} {\bf 24} (2014), 1913-1980.


\bibitem{JP2008} M. Junge, J. Parcet,
Operator space embedding of Schatten $p$-classes into von Neumann algebra preduals,
{\it Geom. Funct. Anal.} {\bf 18} (2008), 522-551.

\bibitem{JP2014} M. Junge, M. Perrin,
Theory of $H_p$-spaces for continuous filtrations in von Neumann algebras,
{\it Ast\'{e}risque} {\bf 362} (2014), vi + 134 pages.

\bibitem{JX2003} M. Junge, Q. Xu,
Noncommutative Burkholder/Rosenthal inequalities,
{\it Ann. Probab.} {\bf 31} (2003), 948-995.

\bibitem{JX2005} M. Junge, Q. Xu,
On the best constants in some noncommutative martingale inequalities,
{\it Bull. London Math. Soc.} {\bf 37} (2005), 243-253.

\bibitem{JX2007} M. Junge, Q. Xu,
Noncommutative maximal ergodic inequalities,
{\it J. Amer. Math. Soc.} {\bf 20}(2) (2007), 385-439.

\bibitem{JX2008} M. Junge, Q. Xu,
Noncommutative Burkholder/Rosenthal inequalities II: Applications,
{\it Israel J. Math.} {\bf 167} (2008), 227-282.

\bibitem{JX2010} M. Junge, Q. Xu,
Representation of certain homogeneous Hilbertian operator spaces and applications,
{\it Invent. Math.} {\bf 175} (2010), 75-118.

\bibitem{KS2008} N. J. Kalton, F. A. Sukochev,
Symmetric norms and spaces of operators,
{\it J. Reine Angew. Math.} {\bf 621} (2008), 81-121.

\bibitem{Lab2013} L. E. Labuschagne,
A crossed product approach to Orlicz spaces,
{\it Proc. Lond. Math. Soc.} {\bf 107} (5) (2013), 965-1003.


\bibitem{LMX2010} C. Le Merdy, Q. Xu,
Maximal theorems and square functions for analytic operators on $L^p$-spaces,
{\it J. London Math. Soc.} {\bf 86}(2) (2012), 343-365.

\bibitem{LMX2011} C. Le Merdy, Q. Xu,
Strong $q$-variation inequalities for analytic semigroups,
{\it Ann .Inst. Fourier.} {\bf 62} (2012), 2069-2097.

\bibitem{LT1979} J. Lindenstrauss, L. Tzafriri,
{\it Classical Banach Space \rm{II},}
Berlin: Springer-Verlag, 1979.



\bibitem{LPP1991} F. Lust-Piquard, G. Pisier.
Noncommutative Khintchine and Paley inequalities.
{\it Arkiv f\"{o}r Mat.} {\bf 29} (1991), 241-260.

\bibitem{LPX2007} F. Lust-Piquard, Q. Xu.
The little Grothendieck theorem and Khintchine inequalities for symmetric spaces of measurable operators.
{\it J. Funct. Anal.} {\bf 244} (2007), 488-503.

\bibitem{M1985} L. Maligranda,
Indices and interpolation,
{\it Dissert. Math.} {\bf 234}, Polska Akademia Nauk, Inst. Mat., 1985.

\bibitem{M1989} L. Maligranda,
{\it Orlicz spaces and interpolation,}
Seminars in Mathematics, Departamento de Matem\'{a}tica, Universidade Estadual de Campinas, Brasil, 1989.


\bibitem{Mei2007} T. Mei,
Operator valued Hardy spaces,
{\it Mem Amer. Math. Soc.} {\bf 881} (2007), 1-64.





\bibitem{PR2006} J. Parcet, N. Randrianantoanina,
Gundy's decomposition for noncommutative martingales and applications,
{\it Proc. London Math. Soc.} {\bf 93}(3) (2006), 227-252.

\bibitem{Perrin2009} M. Perrin,
A noncommutative Davis' decomposition for martingales,
{\it J. London Math. Soc.} (2){\bf 80}(3) (2009), 627-648.


\bibitem{Pisier1998} G. Pisier,
Non-commutative vector valued $L_p$-spaces and completely $p$-summing maps,
{\it Ast\'{e}risque} {\bf 247} (1998), v + 131.

\bibitem{Pisier2003} G. Pisier,
{\it Introduction to Operator Space Theory,}
Cambridge University Press, Cambridge, 2003.

\bibitem{PS2002} G. Pisier, D. Shlyakhtenko,
Grothendieck's theorem for operator spaces,
{\it Invent. Math.} {\bf 150} (2002), 185-217.

\bibitem {PX1997} G. Pisier, Q. Xu,
Non-commutative martingale inequalities,
{\it Commun. Math. Phys.} {\bf 189} (1997), 667-698.

\bibitem{PX2003} G. Pisier, Q. Xu,
Noncommutative $L^p$-spaces,
in {\it Handbook of the Geometry of Banach Spaces, vol.2}: 1459-1517, 2003.

\bibitem{Rand2002} N. Randrianantoanina,
Non-commutative martingale transforms,
{\it J. Funct. Anal.} {\bf 194} (2002), 181-212.


\bibitem{Rand2005} N. Randrianantoanina.
A weak-type inequality for non-commutative martingales and applications.
{\it Proc. London Math. Soc.} {\bf 91}(3) (2005), 509-544.

\bibitem{Rand2007} N. Randrianantoanina,
Conditional square functions for noncommutative martingales,
{\it Ann. Proba.} {\bf 35} (2007), 1039-1070.

\bibitem{RandWu2015a} N. Randrianantoanina, L. Wu,
Martingale inequalities in noncommutative symmetric spaces,
{\it J. Funct. Anal.} {\bf 269} (2015), 2222-2253.

\bibitem{RandWu2015b} N. Randrianantoanina, L. Wu,
Noncommutative Burkholder/Rosenthal inequalities associated with convex functions,
arXiv: 1506.04134, 2015.

\bibitem{RX2014} E. Ricard, Q. Xu,
A noncommutative martingale convexity inequality,
arXiv: 1405.0431.


\bibitem{VDN1992} D. V. Voiculescu, K. J. Dykema, A. Nica,
{\it Free Random Variables,}
Amer. Math. Soc., Providence, RI, 1992.

\bibitem{XXY2015} X. Xiong, Q. Xu, Z. Yin,
Sobolev, Besov and Triebel-Lizorkin spaces on quantum tori,
arXiv: 1507.01789.

\bibitem{Xu1991} Q. Xu,
Analytic functions with values in lattices and symmetric spaces of measurable operators,
{\it Math. Proc. Cambridge Phil. Soc.} {\bf 109} (1991), 541-563.

\bibitem{Xu2006} Q. Xu,
Operator space Grothendieck inequalities for noncommutative $L_p$-spaces,
{\it Duke Math. J.} {\bf 131} (2006), 525-574.

\bibitem{Xu2007} Q. Xu,
Noncommutative $L_p$-Spaces and Martingale Inequalities,
Book manuscript, 2007.

\end{thebibliography}
\end{document}